
\magnification 1200
\input amstex
\documentstyle{amsppt}
\NoBlackBoxes
\NoRunningHeads

\hsize = 6.6 truein

\define\ggm{{G}/\Gamma}
\define\qu{quasiunipotent}
\define\ve{\varepsilon}
\define\hd{Hausdorff dimension}
\define\hs{homogeneous space}
\define\df{\overset\text{def}\to=}
\define\br{\Bbb R}
\define\bc{\Bbb C}

\define\bn{\Bbb N}
\define\bz{\Bbb Z}
\define\bq{\Bbb Q}
\define\ba{badly approximable}
\define\wa{well approximable}
\define\di{Diophantine}
\define\da{Diophantine approximation}

\define\ds{dynamical system}
\define\ca{\Cal A}

\define\h{\goth h}
\define\g{\goth g}
\define\pt{\Phi_t}
\define\un#1#2{\underset\text{#1}\to#2}

\define\edt{e^{\chi t}}
\define\ay{\bold A_\infty}

\define\td{tesselation domain}
\define\tn{tesselation}

\define\ehs{expanding horospherical subgroup}
\redefine\aa{\langle A,\va\rangle}

\define\mr{M_{m,n}(\br)}
\define\mtr{\tilde M_{m,n}(\br)}
\define\amr{$A\in M_{m,n}(\br)$}
\define\amtr{$\aa\in \tilde M_{m,n}(\br)$}
\define\bamn{\Cal B \Cal A_{m,n}}
\define\batmn{\widetilde{\Cal B \Cal A}_{m,n}}

\define\la{L_{A}}
\define\lta{\tilde L_{A,\va}}

\define\va{\bold b}
\define\vx{\bold x}
\define\vc{\bold c}

\define\vv{\bold v}
\define\vu{\bold u}

\define\vq{\bold q}
\define\vp{\bold p}
\define\vw{\bold w}
\define\mn{{m+n}}
\define\nz{\smallsetminus \{0\}}

\redefine\t{^{\sssize T}}


\topmatter
\title Badly approximable systems of affine forms\endtitle

\author { Dmitry Kleinbock} \\ 
  { \rm 
   Rutgers University} 
\endauthor 

    \address{ Dmitry Y. Kleinbock, Department of Mathematics, Rutgers
University, New Brunswick, NJ 08903}
  \endaddress

\email kleinboc\@math.rutgers.edu \endemail

 \thanks Supported in part by NSF
Grants DMS-9304580 and DMS-9704489. \endthanks

\subjclass Primary 11J20, 11J83; Secondary 57S25 \endsubjclass

\abstract We prove an inhomogeneous analogue of W.$\,$M.$\,$Schmidt's
(1969) 
theorem on the \hd\ of the set of \ba\ systems of linear forms. The proof
is based on ideas and methods from the theory of \ds s, in particular,
on abundance of bounded orbits of mixing flows on \hs s of Lie groups. 
\endabstract


\endtopmatter

\document

\heading{\S 1. Introduction}\endheading

\subhead{1.1}\endsubhead For $m,n\in\bn$, we will denote by $\mr$ the
space or real matrices with  $m$ rows and $n$ columns.  
Vectors will be named by  lowercase boldface letters, such as $\vx =
(x_1,\dots,x_k)\t$. 
 $0$ will mean zero vector in any dimension, as well as zero matrix of any size. The norm  $\|\cdot\|$ on $\br^k$ will be always given by $\|\vx\| = \max_{1\le i \le k}|x_i|$. 
 
 All distances (diameters of sets) in various metric spaces will be denoted by ``dist" (``diam"), and $B(x,{r})$ 
will stand for the open ball of radius ${r}$ centered at $x$. 
To avoid confusion, we will sometimes put subscripts indicating the underlying metric space. If the metric space is a group and $e$ is its identity element, we will  write $B({r})$ instead of $B(e,{r})$.
The \hd\ of a subset $Y$ of a metric space $X$ will be denoted by dim$(Y)$, and we will say that $Y$ is {\it thick\/} (in $X$) if   for any nonempty
open subset $W$ of $X$,  
$
\text{dim}(W\cap
Y) = \text{dim}(W)
$ 
(i.e.~$Y$  has full \hd\ at any point of $X$).

A system  of $m$ linear forms in $n$ variables  given by  
\amr\ is called {\it \ba\/} if 
there exists a constant $c>0$ such that for every ${\vp\in\bz^m}$ and
all but finitely many ${\vq\in\bz^n}$ the product 
$\|A\vq + \vp\|^m\|\vq\|^n$ is greater than  $c$;
equivalently, if
$$
c_A \df \liminf_{\vp\in\bz^m,\,\vq\in\bz^n,\,\vq\to\infty}\|A\vq + \vp\|^m\|\vq\|^n > 0\,.
$$ 
Denote by $\bamn$ the set of \ba\ \amr. It has been known since 1920s that $\bamn$ is infinite  (Perron 1921)
and  of zero Lebesgue measure in $\mr$ (Khintchine 1926), and that
$\Cal B \Cal A_{1,1}$ is thick in $\br$ (Jarnik 1929; the latter
result was obtained using continued fractions). In 1969
W.$\,$M.$\,$Schmidt \cite{S3} used the technique of
$(\alpha,\beta)$-games  to show 
that the set
$\bamn$  is thick in $\mr$.

\subhead{1.2}\endsubhead The subject of the present paper is an inhomogeneous analogue of the above notion. By an {\it affine form\/} we will mean a linear form plus a real number. A system  of $m$ affine forms in $n$ variables  will be then given by  
a pair $\aa$, where \amr\ and $\va\in\br^m$. We will denote by $\mtr$ the direct product of $\mr$ and $\br^m$.
Now say that a system  of affine forms  given by  
\amtr\ is {\it \ba\/} if 
$$
\tilde c_{A,\va} \df \liminf_{\vp\in\bz^m,\,\vq\in\bz^n,\,
\vq\to\infty}\|A\vq + \va + \vp\|^m\|\vq\|^n > 0\,,\tag 1.1
$$ 
and {\it \wa\/} otherwise. We will denote by $\batmn$ the set of \ba\ \amtr. 

Before going further, let us consider several trivial examples of \ba\
systems of affine forms. 

\example{1.3. Example} For comparison let us start with the
homogeneous case. Suppose that $A\vq_0  \in\bz^m$ for some 
$\vq_0\in\bz^n\nz$. Then clearly there exist infinitely many
$\vq\in\bz^n$ (integral multiples of $\vq_0$)  for which
$A\vq\in\bz^m$, hence such $A$ is \wa. 
On the other hand, the assumption 
$$
A\vq_0   + \va + \vp_0 = 0\tag 1.2
$$
does not in general guarantee the existence of any other $\vq\in\bz^n$
with $A\vq + \va \in\bz^m$, and, in view of the definition above, just
one integral solution is not enough for $\aa$ to be \wa. 
We will say
that  \amtr\ is {\it rational\/} if
(1.2) holds for some $\vp_0\in\bz^m$ and
$\vq_0\in\bz^n$, and {\it irrational\/} otherwise.

Because of the aforementioned difference of the homogeneous and
inhomogeneous cases, rational systems of forms will have to be treated
separately. In fact, as mentioned in \cite{C, Chapter III, \S 1},
(1.2) allows one to reduce the
study of a rational system \amtr\ to that of $A$. Indeed, for all
$\vq\ne\vq_0$ one can write
$$
\|A\vq + \va + \vp\|^m\|\vq\|^n = \|A(\vq -\vq_0)  + \vp -
\vp_0\|^m\|\vq\|^n = \|A(\vq -\vq_0)  + \vp -
\vp_0\|^m\|\vq - \vq_0\|^n \frac{\|\vq\|^n}{\|\vq - \vq_0\|^n}\,,
$$
which shows that for rational \amtr\ one has $\tilde c_{A,\va} = c_A$;
in particular, $\aa$ is \ba\ iff $A$ is.  
 \endexample

\example{1.4. Example} Another class of examples is given by

\proclaim{Kronecker's Theorem \rm (see  \cite{C, Chapter III,
Theorem IV})}  For \amtr, the following are equivalent:

{\rm (i)} there exists $\ve > 0$ such that for any $\vp\in\bz^m$ and
$\vq\in\bz^n$ one has $\|A\vq + \va + \vp\| \ge \ve$;

{\rm (ii)} there exists $\vu\in\bz^m$ such that $A\t\vu\in\bz^n$ but $\va\t\vu$ is not an integer.
\endproclaim

The above equivalence is straightforward in the $m = n = 1$ case: if
$a = \frac kl \in\bq$, then $\|aq + b + p\| \ge
\text{dist}(b,\frac1l\bz)$, and $a\notin \bq$ implies that $\{(aq + b) \mod
1\}$ is dense in $[0,1]$. In general, it is easy to construct numerous
examples of systems $\aa$
satisfying (ii),  and for such systems one  clearly has $\tilde
c_{A,\va} = +\infty$ in 
view of (i).  Here one notices another difference from the homogeneous case: in
view of Dirichlet's Theorem, one has $c_A < 1$ for any
\amr. \endexample 



\subhead{1.5}\endsubhead It follows from the inhomogeneous version of
Khintchine-Groshev Theorem 
(see \cite{C, Chapter VII, Theorem II}) 
that the set $\batmn$  has Lebesgue measure zero. (See \S 5 for a
stronger statement and other extensions.) A natural
problem to consider is to measure the magnitude of this set  in terms
of the \hd.  One can easily see that the systems of forms $\aa$ which are \ba\
by virtue of 
the two previous examples (that is, either are rational with
$A\in\bamn$ or satisfy the assumption (ii) of Kronecker's Theorem)  
belong to a countable union of proper submanifolds of $\mtr$ and,
consequently, form a set of positive Hausdorff
codimension. Nevertheless, the following is 
true and constitutes the main result of the paper: 

\proclaim{Theorem}  The set  $\batmn$  is thick in $\mtr$.
\endproclaim

This theorem will be proved using results and methods of the paper \cite{KM1}. More precisely, we will derive Theorem 1.5 from Theorem 1.6 below. Before stating the latter, let us introduce some notation and terminology from the theory of Lie groups and homogeneous spaces.

\subhead{1.6}\endsubhead Let $G$ be a connected Lie group, $\g$ its Lie algebra.  Any $X\in\g$ gives rise to a one-parameter semigroup $F = \{\exp(tX)\mid t\ge 0\}$, where $\exp$ stands for the exponential map from $\g$ to $G$. We will be interested in the left action of $F$ on \hs s $\Omega\df\ggm$, where $\Gamma$ is a discrete subgroup of $G$. 

Many properties of the above action can be understood by looking at the adjoint action of $X$ on $\g$.  For $\lambda\in\bc$, we denote by $\g_\lambda(X)$ the {\it generalized eigenspace\/} of $\text{ad}\,X$ corresponding to $\lambda$, i.e.~the subspace of the complexification $\g_{\sssize\bc}$ of $\g$ defined by 
$$
\g_\lambda(X) = \{Y\in \g_{\sssize {\Bbb C}}\mid (\text{ad}\,X - \lambda I)^jY =
0\text{ for some }j\in\bn\}\,.
$$
We will say that $X$ is {\it semisimple\/} if $\g_{\sssize\bc}$ is spanned by eigenvectors of $\text{ad}\,X$.
Further, we will define the $X$- (or $F$-) {\it expanding horospherical subgroup\/} of $G$ 
as follows: $H = \exp \h$, where $\h$ is the subalgebra of 
 $\g$  with complexification
$
\h_{\sssize {\Bbb C}} = \underset{\text{Re}\,\lambda > 0}\to\oplus\,\g_\lambda(X)$. 

Say that a discrete subgroup $\Gamma$ of $G$ is a {\it lattice\/} if the quotient space $\Omega = \ggm$ has finite volume with respect to a $G$-invariant measure. Note that $\Omega$ may or may not be compact. Any group admitting a lattice is unimodular; we will choose a Haar measure $\mu$ on $G$ and the corresponding Haar measure $\bar \mu$ on $\Omega$ so that $\bar\mu(\Omega) = 1$. The $F$-action on $\Omega$ is said to be {\it mixing\/} if 
$$
\lim_{t\to+\infty}\bar\mu\big(\exp(tX)W \cap K\big) = \bar\mu(W)\bar\mu(K)\tag 1.3
$$
for any two measurable subsets $K$, $W$ of $\Omega$. 

The last piece of notation comes from the papers \cite{K2, K3}. Consider the one-point compactification $\Omega^*\df\Omega\cup\{\infty\}$ of $\Omega$, topologized so that the complements
to all compact sets constitute the basis of neighborhoods of
$\infty$. We will use the notation $Z^*\df Z\cup\{\infty\}$ for any subset $Z$  of $\Omega$. Now for a subset $W$ of $\Omega^*$ and a subset $F$ of $G$ define
$E(F,W)$ to be the set of points of $\Omega$ with $F$-orbits escaping $W$, that is
$$
E(F,W) \df \{x \in \Omega\mid \overline{Fx}\cap
W = \varnothing \}\,,
$$  
with the closure taken in the topology of $\Omega^*$. In particular, if $Z$  is a subset of $\Omega$, $E(F,Z^*)$ stands for the set of $x\in\Omega$ such that orbits $Fx$ are bounded and stay away from $Z$. 

We are now ready to state 

\proclaim{Theorem}  Let $G$ be a real Lie group,  $\Gamma$ a lattice in $G$, $X$ a semisimple element of the Lie algebra of $G$, $H$ the $X$-expanding horospherical subgroup of $G$ 
such that the action of $F = \{\exp(tX)\mid t\ge 0\}$ on $\Omega = \ggm$ is mixing. Then for any closed $F$-invariant null
subset $Z$ of $\Omega$ and any $x\in\Omega$, the set
$
\{h \in H\mid hx \in E(F,Z^*)\}$ is thick in $H$. 
\endproclaim

\subhead{1.7}\endsubhead The reduction of Theorem 1.5 to the above
theorem is described in \S 4 and is based on a version of
S.$\,$G.$\,$Dani's (see 
\cite{D1} or \cite{K3, \S 2.5}) correspondence between  \ba\ systems
of linear forms and certain orbits of lattices in the Euclidean space
$\br^\mn$. More precisely, let $G_0 = SL_\mn(\br)$ and let $G$ be the
group $\text{\rm Aff}(\br^\mn)$ of measure-preserving affine
transformations of $\br^\mn$, i.e. the semidirect product of $G_0$ and
$\br^\mn$. Further, put  $\Gamma_0 = SL_\mn(\bz)$ and let $\Gamma =
\text{\rm Aff}(\bz^\mn) \df \Gamma_0\ltimes\bz^\mn$ be the subgroup of
$G$ leaving the standard lattice $\bz^\mn$ invariant. It is easy to
check that $\Gamma$  is a non-cocompact lattice in
$G$, and that $\Omega \df \ggm$ can be identified with the space of
{\it free unimodular lattices\/} in $\br^\mn$, i.e.  
$$
\Omega\cong\{\Lambda +\vw\mid \Lambda\text{ a lattice in $\br^\mn$ of covolume $1$},\ \vw\in\br^\mn\}\,.
$$ 
One can show that the quotient topology on $\Omega$
coincides with the natural topology on the space of lattices, so that
two lattices are close to each other if their generating elements
are. We will also consider $\Omega_0 \df G_0/\Gamma_0$ and identify it
with  the subset of $\Omega$ consisting of
``true'' lattices, i.e.~those containing $0\in\br^\mn$. 

We will write elements of $G$ in the form $\langle L,\vw\rangle $, where $L\in G_0$ and $\vw\in\br^\mn$, so that $\langle L,\vw\rangle $ sends $\vx\in\br^\mn$ 
to $L\vx + \vw$. If $\vw = 0$, we will simply write $L$ instead of
$\langle L,0\rangle $; the same convention will apply to elements of
the Lie algebra $\g$ of $G$. We will fix an element $X$ of
$\goth{sl}_\mn(\br)\subset \g$ of the form  
$$
X = \text{diag}(\underbrace{\tfrac1m,\dots,\tfrac1m}_{\text{$m$ times}}, \underbrace{-\tfrac1n,\dots,-\tfrac1n}_{\text{$n$ times}})\,,
$$
and let $F = \{\exp(tX)\mid t\ge 0\}$. 

Recall that in the standard version of Dani's correspondence, to a system of linear forms given by \amr\ one associates a lattice $\la\bz^\mn$ in $\br^\mn$, where $\la$ stands for $\left(\matrix
I_m & A  \\
0 & I_n
\endmatrix \right)$. It is proved in \cite{D1} that $A \in \bamn$ if
and only if the trajectory $F\la\bz^\mn\subset \Omega_0$ is bounded
(see \S 4.3 for more details).
From this and the aforementioned 1969 result of Schmidt, Dani
deduced 
in \cite{D1} that   the set of lattices in $\Omega_0$ with bounded
$F$-trajectories 
is thick. 

It was suggested by Dani  and then conjectured by G.$\,$A.$\,$Margulis
\cite{Ma, Conjecture (A)} that the abundance of bounded orbits is a
general feature of non\qu\ (see \S2.3) flows on \hs s of Lie
groups. This conjecture was settled in \cite{KM1}, thus giving an
alternative (dynamical) proof  of Dani's result, and hence of the
thickness of the set $\bamn$.

Our goal in this paper is to play the same game in the inhomogeneous
setting. Given \amtr\ one can consider a free lattice in $\br^\mn$ of
the form $\lta\bz^\mn$, where $\lta$ is the element of $G$ given by
$\langle \la,(\va,0)\t\rangle $.  
In \S 4 we will prove the following

\proclaim{Theorem}  Let $G$, $\Gamma$, $\Omega$, $\Omega_0$ and $F$ be
as above. Assume that   $\lta\bz^\mn$ 
belongs to \linebreak $E\big(F,(\Omega_0)^*\big)$. Then   a system of affine  forms given by $\aa$
is 
\ba\.  \endproclaim 

To see that this theorem provides a link from Theorem 1.6 to Theorem 1.5, it remains to observe that $F$-action on $\Omega$ is mixing (all the necessary facts related to mixing of actions on \hs s are collected in \S 2), and that $\{\lta\mid\text{\amtr}\}$ is the $F$-expanding horospherical subgroup of $G$. 
In fact, Theorem 1.7 is obtained as a corollary from a necessary and sufficient condition for an irrational system of affine forms to be badly approximable, an inhomogeneous analogue of Dani's correspondence
 \cite{D1, Theorem 2.20} and the main result of \S 4 of the present paper. 

For the sake of making this paper self-contained, in \S 3 we present
the complete proof of Theorem 1.6, which is
basically a simplified version of the argument from \cite{KM1}. The
last section of the paper is devoted to several concluding remarks and
open questions.   

\heading{\S 2. Mixing and the \ehs}\endheading

\subhead{2.1}\endsubhead Throughout the next two sections, we let $G$
be a connected Lie group, $\g$ its Lie algebra, $X$ an element of
$\g$, $g_t = \exp(tX)$,  $F = \{g_t\mid t\ge 0\}$, $\Gamma$ a lattice
in $G$ and $\Omega = \ggm$. For  $x\in\Omega$, denote by $\pi_x$  the
quotient map $G\to\Omega$, $g\to gx$. The following restatement of the
definition of mixing of the $F$-action on $\Omega$ is straightforward:

\proclaim{Proposition \rm (cf.~\cite{KM1, Theorem 2.1.2})} The action
of $F$ on $\Omega$ is mixing iff for any  compact $Q\subset \Omega$,
any measurable $K\subset \Omega$ and any measurable $U\subset G$ such
that $\pi_x$ is injective on $U$ for all $x\in Q$, one has 
$$
\forall\,\ve>0\ \exists\, T>0\text{ such that }
\left| \bar\mu\big(g_tUx \cap K\big) - \mu(U)\bar\mu(K)
\right|\le\ve \text{ for all } t\ge T\text{ and }x\in Q\,.\tag 2.1
$$ 
\endproclaim

\demo{Proof} To get (1.3) from (2.1), take any $x\in \Omega$, put $Q = \{x\}$ and take $U$ to be any one-to-one $\pi_x$-preimage of $W$. For the converse, one considers the family of sets $W_x = Ux$ and observes that the difference  $\bar\mu\big(g_tW_x \cap K\big) - \bar\mu(W_x)\bar\mu(K)$ goes to zero uniformly when $x$ belongs to a compact subset of $\Omega$. \qed
\enddemo

\subhead{2.2}\endsubhead If $G$ is a  connected semisimple  Lie group
without compact factors 
and $\Gamma$ an irreducible lattice in $G$, one has C.$\,$Moore's
\cite{Mo} criterion for mixing of one-parameter subgroups of $G$:
an $F$-action on $\Omega$ is mixing iff $F$ is not relatively compact in
$G$. Since the group $\text{\rm Aff}(\br^\mn)$  is not semisimple, we
need a reduction to the 
semisimple case based on the work \cite{BM} of Brezin and Moore.
Following \cite{Ma}, say that a \hs\ $G/\Delta$ is a {\it quotient\/} of
$\Omega$ if $\Delta$ is a closed subgroup of $G$ containing $\Gamma$.  If $\Delta$
contains a closed normal subgroup $L$ of $G$ such that $G/L$ is
semisimple (resp.~Euclidean\footnote{A connected solvable Lie group is called {\it
Euclidean\/} if it is locally isomorphic to an extension of a vector
group by a compact Abelian Lie group. }) then the
quotient $G/\Delta$ is called 
{\it semisimple\/} (resp.~{\it Euclidean\/}). It is easy to show that
the maximal semisimple (resp.~Euclidean) quotient
of $\Omega$ exists (by the latter we mean the quotient  $G/\Delta$ of
$\Omega$ such that any other  semisimple (resp.~Euclidean)
quotient of $\Omega$ is a quotient of $G/\Delta$). 

The following
proposition is a
combination of Theorems 6 and 9 from \cite{Ma}:
 
\proclaim{Proposition} Suppose
that 

{\rm (i)} there are no nontrivial  Euclidean quotients of $\Omega$, and 

{\rm (ii)}  $F$ acts ergodically on the maximal semisimple quotient of
$\Omega$.

 Then the action of $F$ on $\Omega$ is mixing. \endproclaim

\subhead{2.3}\endsubhead Choose a Euclidean structure on $\g =
\text{Lie}(G)$, inducing a right-invariant Riemannian metric on $G$
and a corresponding  Riemannian metric on $\Omega$. We will fix a
positive $\sigma_0$ such that 
$$
\aligned
\text{the }&\text{restriction of $\exp:\g\to G$ to  $B_\g(4\sigma_0)$ is
one-to-one}\\
&\text{and distorts distances by at most a factor of $2$}\,.
\endaligned\tag 2.2
$$
Denote by $\tilde \h$ the subalgebra of 
 $\g$  with complexification
$
\tilde \h_{\sssize {\Bbb C}} = \underset{\text{Re}\,\lambda > 0}\to\oplus\g_\lambda(X)$, and put $\tilde H = \exp \tilde \h$. Clearly  $\g = \h \oplus \tilde \h$, which implies that  the multiplication map $\tilde H \times H\to G$ is one-to-one in a neighborhood of the identity in $\tilde H \times H$.

We now assume that  the $X$-\ehs\ $H$ of $G$ is nontrivial (in the
terminology of \cite{Ma}, $F$ is not {\it \qu\/}). Denote by $\chi$
the trace of $\text{ad}\,X|_\h$ and by $\lambda$ the real part of an
eigenvalue of  $\text{ad}\,X|_\h$ with the smallest real part. Denote
also by $\pt$ the inner automorphism $g\to g_tgg_{-t}$ of $G$. Since
$\text{Ad}\,g_t = e^{\text{ad}\,tX}$ is the differential of $\pt$ at
the identity,   the  Jacobian of $\pt|_{\sssize H}$ is equal to
$e^{\chi t}$, and local metric properties of $\pt$ are determined by
eigenvalues of $\text{ad}\,tX$.   In particular, 
the following is true: 

\proclaim{Lemma} For all $t > 0$ one has

\rm{(a)} ${\text{\rm dist}\big(\pt^{-1}(g),\pt^{-1}(h)\big)}
\le 4 e^{-\lambda t}{\text{\rm dist}(g,h)}\text{ for all }g,h\in  B_{\sssize H}(\sigma_0)$;

\rm{(b)} if $X$ is semisimple, ${\text{\rm dist}\big(\pt(g),\pt(h)\big)}
\le 4 \cdot {\text{\rm dist}(g,h)}\text{ for all }g,h\in  B_{\sssize \tilde H}(\sigma_0)$.
\endproclaim

In other words, $\pt$ acts as an expanding map of $H$ and as a non-expanding map of $\tilde H$.

\subhead{2.4}\endsubhead We now turn to a crucial application of mixing of $F$-action on $\Omega$. Choose a Haar measure $\nu$ on $H$. Roughly speaking, our goal is to replace a subset $U$ of $G$ in the formula (2.1) with a $\nu$-measurable subset $V$ of $H$. 

\proclaim{Proposition} Let $V$ be a bounded measurable subset of $H$, $K$ a bounded measurable subset of $\Omega$
with $\bar\mu(\partial K) = 0$,   
$Q$ a compact subset of $\Omega$. Assume that $X$ is semisimple and that the $F$-action on $\Omega$ is mixing. Then for any $\ve>0$  there exists
$T_1 = T_1(V,K,Q,\ve)>0$ such that
$$
t\ge T_1\Rightarrow\forall\, x\in Q\quad
\nu\big(\{h\in V\mid g_thx\in K\}\big) \ge
\nu(V)\bar\mu(K) - \ve\,. \tag 2.3
$$
\endproclaim

\demo{Proof} Since $V$ is bounded,  $Q$ is compact, and $\Gamma$ is discrete,
$V$ can be decomposed as a disjoint union of subsets $V_j$ of $H$ 
 with $\pi_x$ injective on some neighborhood of $V_j$ for
all $x\in Q$ and for each $j$. Hence one can without loss of
generality assume that the maps
$\pi_x$ are injective on some neighborhood $U'$ of $V$ for all
$x\in Q$. Similarly, one can safely assume that $V\subset B_{\sssize H}(\sigma_0)$ and $\nu(V) \le 1$.

Choose a subset $K'$ of $K$ such that $\bar\mu(K') \ge \bar\mu(K) - \ve/2$ and the distance $\sigma_1$ between $K'$ and $\partial K$ is positive. Then choose $\sigma\le \min(\sigma_0, \sigma_1/4)$. After that, pick a neighborhood $\tilde V\subset B_{\sssize \tilde H}(\sigma)$ of the identity in $\tilde H$
such that  $U\df\tilde VV$ is contained in $U'$. 

Given $x\in\Omega$ and $t > 0$, denote by $V'$ the set $\{h\in V\mid g_thx\in K\}$ that we need to estimate the measure of. 

\proclaim{Claim} The set $Ux\cap g_t^{-1}K'$ is contained in $\tilde VV'x$. \endproclaim

\demo{Proof} For any $h\in H$ and $\tilde h\in\tilde V$,  $g_t\tilde hhx = \pt(\tilde h)g_thx\in B(g_thx, 4\sigma)\subset B(g_thx, \sigma_1)$ by Lemma 2.3(b) and the choice of $\tilde V$. Therefore  $g_thx$ belongs to $k$ whenever $g_t\tilde hhx\in K'$. \qed
\enddemo

Now, using Proposition 2.1, find $T_1> 0$ such that $\left| \bar\mu\big(g_tUx \cap K'\big) - \mu(U)\bar\mu(K')
\right|\le\ve/2$  for all  $t\ge T_1$ and $x\in Q$. In order to pass from $U$ to $V$, 
choose a left Haar measure $\tilde \nu$ on $\tilde H$ such that $\mu$ is the product of $\nu$ and $\tilde \nu$ (cf.~\cite{Bou,
Ch.~VII, \S 9, Proposition 13}). Then for all  $t\ge T_1$ and $x\in Q$ one can  write
$$
\split
&\tilde \nu(\tilde V) \nu(V') = \mu(\tilde VV') = \bar\mu(\tilde VV'x) \un{(Claim)}\ge \bar\mu(Ux\cap g_t^{-1}K') = \bar\mu(g_tUx\cap K') \\
\un{(if $t\ge T_1$)}\ge &\mu(U)\bar\mu(K') - \ve/2 
\ge \tilde \nu(\tilde V) \nu(V)\big(\bar\mu(K) - \ve/2\big) -\ve/2 \ge \tilde \nu(\tilde V) \big(\nu(V)\bar\mu(K)  -\ve\big)\,,
\endsplit
$$
which immediately implies (2.3).  
 \qed
\enddemo

Note that similarly one can estimate $\nu(V')$ from above,
see \cite{KM1, Proposition 2.2.1} for a more general statement. 

It will be convenient 
to denote the $\pt$-image of $V'$ by $V(x,K,t)$,
i.e.~to let 
$$
V(x,K,t) \df  \pt\big(\{h\in V\mid g_thx\in K\}\big) = \{h\in
\pt(V)\mid hg_tx\in 
K\}\,. 
$$
Roughly speaking, Proposition 2.4 says that the relative measure of $
V(x,K,t)$ in $\pt(V)$ is big when $t$ is large enough: indeed, (2.3)
can be rewritten in the form
$$
t\ge T_1\Rightarrow\forall\, x\in Q\quad
\nu\big(V(x,K,t)\big) \ge
\edt\big(\nu(V)\bar\mu(K) - \ve\big)\,. \tag 2.4
$$

%

\heading{\S 3. Proof of Theorem 1.6}\endheading 


\subhead{3.1}\endsubhead  Following \cite{KM1}, say  that an
open subset $V$ of $H$ is a {\it \td} for the right action of $H$ on
itself  relative to a
countable subset $\Lambda$ of $H$ if

(i) $\nu (\partial V) = 0\,,$

(ii) $V\gamma_1 \cap V\gamma_2  = \varnothing$ for different $\gamma_1,
\gamma_2 \in \Lambda\,,$ and

(iii) $X = \bigcup_{\gamma\in\Lambda}\overline{V}\gamma\,$.

The pair $(V,\Lambda)$ will be called a {\it
\tn} of $H$. Note that it follows easily from (ii) and (iii) that for
any measurable subset $A$ of $H$ one has 
$$
\frac{\nu(A)}{\nu( V )} \le \#\{\gamma\in\Lambda\mid 
V\gamma\cap A\ne\varnothing\} 
\le \frac{\nu\big(\{h\in H\mid\text{dist}(h,A)\le \text{diam}(V)\}\big)}{\nu( V
)}\,.\tag 3.1   
$$

Let $k$ stand for the dimension of $H$. We will use a one-parameter family of \tn s of $H$ defined as follows: if 
 $\{X_1, \dots, X_k\}$ is a fixed orthonormal {\it strong Malcev basis\/} of $\h$ (see \cite{CG} or \cite{KM1, \S 3.3} for the definition), we let $I =  \big\{\sum_{j=1}^k
x_jX_j\bigm| |x_j|<1/2\big\}$ be the unit cube  in $\h$, and then
take $V_r=\exp\big(\frac{r}{\sqrt{k}}I\big)$. It was proved in
\cite{KM1} that $V_r$ is a \td\ of $H$; let $\Lambda_r$ be a
corresponding set of translations. It is clear from (2.2) that $
V_{r}$ is contained in $B({r})$ provided $r\le\sigma_0$, where  $\sigma_0$
is as in  \S 2.3.

The main ingredient of the proof of Theorem 1.6 is given by the
following procedure: we look at the expansion $\pt(V_r)$ of the set
$V_r$  by the
automorphism $\pt$, and then consider the translates 
$V_{r}\gamma$ which lie entirely inside $\pt(V_r)$. 
It was shown in \cite{KM1} (see also \cite{K2, Proposition 2.6}) that
when $t$ is large enough, the measure of the union of all such translates
 is approximately equal to
the measure of  $\pt(V_{r})$; in other words, boundary effects are
negligible.  More precisely, the following is what will be needed for
the proof of the main theorem:

\proclaim{Proposition} For any $r\le\sigma_0$  and any $\ve>0$ there exists $T_2 =
T_2(r,\ve)>0$ such that 
$$
t\ge T_2\quad\Rightarrow\quad\#\big\{\gamma\in\Lambda_{r}\mid 
V_{r}\gamma\cap\partial\big(\pt(V_{r})\big) \ne \varnothing\big\}
\le \ve e^{\chi t}\,.
$$
\endproclaim

\demo{Proof} One can write 
$$
\#\big\{\gamma\in\Lambda_{r}\mid 
V_{r}\gamma\cap\partial\big(\pt(V_{r})\big) \ne \varnothing\big\} =
\#\big\{\gamma\in\Lambda_{r}\mid \pt^{-1}(V_{r}\gamma)\cap\partial V_{r} \ne \varnothing
\big\}\,.
$$
Observe that $\big(\pt^{-1}(V_{r}),\pt^{-1}(\Lambda_{r})\big)$ is also a
\tn\ of $H$, and, in view of Lemma 2.3(a), the diameter of
$\pt^{-1}(V_{r})$ 
is at most  $8re^{-\lambda t}$. Therefore,  by (3.1),  the number
in the right hand 
side is not greater than the ratio of the measure of the $8re^{-\lambda
t}$-neighborhood of $\partial V_{r}$ (which, in view of condition (i)
above, tends to zero as 
$t\to\infty$) and $\nu\big(\pt^{-1}(V_{r})\big)
= e^{-\chi
t}\nu(V_r)$. This shows that $\lim_{t\to\infty}
e^{-\chi t}\#\big\{\gamma\in\Lambda_{r}\mid 
V_{r}\gamma\cap\partial\big(\pt(V_{r})\big) \ne
\varnothing\big\} = 0$, hence the proposition.  
 \qed
\enddemo

\subhead{3.2}\endsubhead Suppose  a  subset $K$ of
$\Omega$, a point $x\in\Omega$, $t > 0$ and positive ${r}\le\sigma_0$
are given.  Consider a \tn $(V_r,\Lambda_r)$ of $H$, and recall that
we defined $V_r(x,K,t)$ as the set of all 
elements $h$ in $\pt(V_r)$ for which $h g_tx$ belongs to $K$. Our goal
now is to approximate this set by  the union of translates of
$V_r$. More precisely, let us denote by $\Lambda_r(x,K,t)$ the set of 
translations $\gamma \in\Lambda_r$ such that $V_r\gamma$ lies entirely inside
$V_r(x,K,t)$; in other words, if $V_r\gamma\subset\pt(V_r)$ and 
$V_r\gamma g_tx = g_t\pt^{-1}(V_r\gamma)x$ is contained in $K$. Then
the union  
$$
\bigcup_{\gamma\in\Lambda_{{r}}(x,K,t)} V_r\gamma =
\bigcup_{\gamma\in\Lambda,\, V_r\gamma\subset V_r(x,K,t)} V_r\gamma
$$
 can be 
thought of as a ``\tn\ approximation'' to $V_r(x,K,t)$.
 We can
therefore think of the theorem below as of a ``\tn\ approximation'' to 
 Proposition 2.4.

\proclaim{Theorem} Let $K$ be a subset of $\Omega$ with
$\bar\mu(\partial K) = 0$, $Q$ a compact subset of $\Omega$. Then for
any $\ve > 0$  there exists $r_0 = r_0(K,\ve) \in(0,\sigma_0)$ such that 
for any positive $r\le r_0$ one can find $T_0 = T_0(K,Q,\ve,r)>0$ with
the following property:
$$
 t\ge T_0\Rightarrow\forall\, x\in Q \quad \#\Lambda_{{r}}(x,K,t)\ge
\edt\big(\bar\mu(K) - \ve\big)\,.\tag 3.2
$$
\endproclaim

\demo{Proof} If $\bar\mu(K) = 0$, there is nothing to prove.
Otherwise, pick a compact subset $K'$ of $K$ with $\bar\mu(\partial
K') = 0$, which satisfies $\bar\mu(K')\ge\bar\mu(K)-\ve/3$ and lies at
a positive distance from the complement of $K$. Take $r_0 \le \sigma_0$
such that $V_{r_0} V_{r_0}^{-1}K'\subset K$ (hence $V_r
V_r^{-1}K'\subset K$ for any positive $r\le r_0$). Then for any $t > 0$
and $x\in\Omega$ one has 
$$
V_r\gamma g_tx \subset K 
\ \Leftarrow\  V_r\gamma g_tx \subset V_r V_r^{-1}K' \ \Leftarrow\   \gamma
g_tx\in V_r^{-1}K' \ \Leftarrow\   V_r\gamma
g_tx\cap K' \ne \varnothing \,.
$$
Therefore
$$
\split
\#\Lambda_{{r}}(x,K,t) &=\#\{\gamma\in\Lambda_{{r}}\mid V_r\gamma
\subset \pt(V_r) \ \&\  V_r\gamma g_tx \subset  K \}\\
&\ge \#\{\gamma\in\Lambda_{{r}}\mid V_r\gamma
\subset \pt(V_r) \ \&\  V\gamma
g_tx\cap K' \ne \varnothing\}\\
&\ge \#\{\gamma\in\Lambda_{{r}}\mid V_r\gamma \cap V_r(x,K',t) \ne
\varnothing\} - \#\{\gamma\in\Lambda_{{r}}\mid V_{r}\gamma\cap\partial\big(\pt(V_{r})\big) \ne \varnothing\}\,.
\endsplit
$$

 Now take  
$$
T_0 = \max\, \left(T_1\big(V_r,K',Q,\frac{{\ve}\nu(V)}3\big)\text{ from
Proposition 
2.4, \  $T_2\big(r,\frac\ve3\big)$ from Proposition
3.1}\right)\,.
$$
 Then   
the number of $\gamma\in\Lambda_{{r}}$ for which $V_r\gamma$ has
nonempty intersection with $V_r(x,K',t)$  is, in view of (2.4) and
(3.1), for all $x\in Q$ and 
$t\ge T_0$ not less than 
$$
\frac{\edt\big(\nu(V_r)\bar\mu(K') - \ve\nu(V_r)/3\big)}{\nu(V_r)} =
{\edt\big(\bar\mu(K') - \ve/3\big)} \ge {\edt\big(\bar\mu(K) - 2\ve/3\big)}\,.
$$ 
 On the other hand, the  number of translates $V_r\gamma$ nontrivially
intersecting with $\partial\big(\pt(V_{r})\big)$ is, by
Proposition 
3.1, not greater than $\ve e^{\chi t}/3$, and (3.2) follows.
\qed\enddemo

\subhead{3.3}\endsubhead  We now describe a construction  of a class of sets  
for which there is a natural lower estimate for the \hd.   Let $X$ be a
Riemannian manifold, $\nu$ a 
Borel measure on $X$, 
$A_0$ a compact subset of $X$. Say that a countable
collection $\ca$ of compact subsets of $A_0$ of positive measure $\nu$ is {\it tree-like\/}
relative to $\nu$ if
$\ca$ is the union of finite nonempty subcollections $\ca_j$,
$j = 0,1,\dots$, such
that $\ca_0 = \{A_0\}$ and the following two conditions are satisfied:
$$
\align
\forall\,j\in\Bbb N\quad \forall\, A, B \in \ca_j\quad&\text{either
}A=B\quad\text{or}\quad \nu(A\cap B) = 0\,;\tag TL1 \\
\forall\,j\in\Bbb N\quad \forall\,B \in \ca_j\quad\,&\exists\,A\in
\ca_{j-1} \quad\text{such that}\quad B\subset A\,.\tag TL2
\endalign
$$
Say also that $\ca$ is {\it strongly tree-like\/} if it is tree-like and
in addition
$$
d_j(\ca)\df\sup_{A\in\ca_j}\text{diam}(A)\to 0\quad\text{as}\quad
j\to\infty\,.\tag STL 
$$

Let $\ca$ be a tree-like collection of sets. 
For each $j = 0,1,\dots$, let ${\bold A}_j = \bigcup_{A\in \ca_j}A$. These are
nonempty compact sets, and from (TL2) it follows that
${\bold A}_j\subset{\bold A}_{j-1}$ for any $j\in \Bbb N$. Therefore one can define
the (nonempty) {\it limit set\/} of $\ca$ to be 
$$
\ay = \bigcap_{j = 0}^\infty{\bold A}_j\,.
$$ 

Further, for any subset $B$ of $A_0$ with $\nu(B) > 0$ and
any $j\in \Bbb N$, define the {\it $j$th 
stage density\/} $\delta_j(B,\ca)$ of $B$ in $\ca$ by 
$$
\delta_j(B,\ca) = \frac{\nu({\bold A}_j\cap B)}{\nu(B)}\,,
$$
and 
the {\it
$j$th stage density\/} $\delta_j(\ca)$ of $\ca$ by
$
\delta_j(\ca) = \inf_{B\in\ca_{j-1}}\delta_{j}(B,\ca)\,.
$

The following estimate, based on an application of Frostman's Lemma,
is essentially proved in \cite{Mc} and \cite{U}:

\proclaim{Lemma} Assume that there exists $k>0$
such that 
$$
\liminf_{r\to 0}\frac{\log
{\nu\big(B(x,r)\big)}}{\log {r}}\ge k\tag 3.3
$$
for any $x\in A_0$. 
Then for any strongly tree-like (relative to $\nu$)
collection $\ca$ of subsets of $A_0$,
$$
\text{\rm dim}(\ay)\ge k -
\limsup_{j\to\infty}\frac{\sum_{i=1}^{j}\log{\delta_i(\ca)}}{\log{d_j(\ca)}}\,.
$$
\endproclaim

\subhead{3.4}\endsubhead Now everything is ready for the

\demo{Proof of Theorem 1.6}  Let $x\in\Omega$ and a nonempty  open subset
$V$ of $H$ be given. We need to prove that the \hd\ of the set $
\{h \in V\mid hx \in E(F,Z^*)\}$ is equal to $k = $dim$(H)$. 
Replacing $x$ by $hx$ for some $h\in V$ we can assume that
 $V$ is a neighborhood of identity in $H$.

 Pick a compact set
$K\subset\Omega \smallsetminus Z$ with $\bar\mu(\partial K) = 0$, 
and choose arbitrary $\ve > 0$, $\ve < \mu(K)$.  
Then, using Theorem 3.2, find $r \le r_0(K,\ve)$ such that the
corresponding \td\ $\overline V_{r}$ is
contained in $V$, and then take $t \ge
\max\big(T_0(K, K\cup \{x\},\ve,r),1/\ve)$. We claim that 
$$
\text{dim}\big(\{h \in \overline{V_r}\mid hx \in E(F,Z^*)\}\big) \ge  k -
\frac{\log\big(\frac 
1{\bar\mu(K)-\ve}\big)}{\lambda t - \log 4}\,.\tag 3.4
$$
Since $\overline{V_r}\subset V$, $\ve$ is arbitrary small and $t$ is
greater than $1/\ve$, it follows from (3.4) that $
\text{dim}\big(\{h \in V\mid hx \in E(F,Z^*)\}\big)$ is equal to $k$.

To demonstrate (3.4), for all $y\in K\cup \{x\}$ let us define strongly
tree-like (relative to the Haar 
measure $\nu$ on $H$) collections $\ca(y)$ inductively as 
follows. We first let $\ca_0(y) = \{\overline{V_r}\}$ for all $y$, 
then define 
$$
\ca_{1}(y) =
\{\Phi_{-t}(\overline{V_r}\gamma) \mid \gamma \in
\Lambda_r(y,K,t)\}\,.\tag 3.5
$$
More generally, if $\ca_i(y)$ is defined for all $y\in K\cup \{x\}$
and $i < j$, 
we let
$$ 
\ca_{j}(y) =
\{\Phi_{-t}(A\gamma) \mid \gamma \in
\Lambda_r(y,K,t),\ A\in\ca_{j-1}(\gamma g_{t}y)\}\,.\tag 3.6
$$
By definition, $\gamma \in
\Lambda_r(y,K,t)$
 implies that 
$\gamma
g_{t}y\in K$; therefore $\ca_{j-1}(\gamma
g_{t}y)$ in (3.6) is defined and the inductive procedure goes through.
The properties (TL1) and (TL2) follow readily from the construction and
$V_r$ being a \td. Also, by Lemma 2.3(a), the diameter of
$\Phi_{-t}(A\gamma)$ is not greater than  $4 e^{-\lambda 
t}$diam$(A)$, which implies that 
$d_j\big(\ca(y)\big)$ is  for all $j\in \bn$ and $y\in K$ not greater
than $2r\cdot (4 e^{-\lambda
t})^j$, and therefore (STL) is satisfied. 

Let us now show by induction that the $j$th stage density
$\delta_j\big(\ca(y)\big)$ of $\ca(y)$ is for all  $y\in K\cup \{x\}$ and $j\in
\bn$ bounded 
from below by $\bar\mu(K)-\ve$.
 Indeed, by definition
$$
\align
\delta_1\left(\overline{V_{r}},\ca(y)\right) &=
\frac{\nu\big({\bold A}_1(y)\big)}{\nu\left(\overline V_{r}\right)} 
\underset{\text{(by (3.5))}}\to=
\frac{\nu\big(\bigcup_{\gamma\in\Lambda_{{r}}(y,K,t)}
\Phi_{t}^{-1}({V_{r}}\gamma)\big)}{\nu(V_{r})}\\ 
&=
\quad e^{-\chi t}\#\Lambda_{{r}}(y,K,t)
 \underset{\text{(by (3.2))}}\to\ge \bar\mu(K)-\ve\,.  
\endalign
$$
On the other hand, if $j\ge 2$ and $B\in \ca_{j-1}(y)$ is of the form $\Phi_{t}^{-1}(A\gamma)$ for  $A\in\ca_{j-2}(\gamma
g_{t}y)$, the formula (3.6) gives
$$
\align
\delta_j\big(B,\ca(y)\big) = \frac{\nu\big(B\cap {\bold A}_j(y)\big)}{\nu(B)} &= \frac{\nu\left(\Phi_{t}^{-1}\big(B\cap {\bold A}_j(y)\big)\right)}{\nu\big(\Phi_{t}^{-1}(B)\big)} \\
= \frac{\nu\big(A\gamma\cap {\bold A}_{j-1}(\gamma g_{t}y)\gamma\big)}{\nu(A\gamma)} &= 
\frac{\nu\big(A\cap {\bold A}_{j-1}(\gamma g_{t}y)\big)}{\nu(A)}  = \delta_{j-1}\big(A,\ca(\gamma g_{t}y)\big)\,,
\endalign
$$
and induction applies. Finally, the measure $\nu$ clearly
satisfies (3.3) 
with  $k=\text{dim}(H)$, 
and an application of Lemma
3.3 yields that for all  $y\in K\cup \{x\}$ one has
$$
\text{\rm
dim}\big(\ay(y)\big)
\ge k - \limsup_{j\to\infty}\frac{j\log\big(\bar\mu(K)-\ve\big)}{
\log\big(2r\cdot (4 e^{-\lambda
t})^j\big)}\,,
$$
which is exactly the right hand side of (3.4). 

To finish the proof it remains to show that $\ay(x)x$ is a subset of
$E(F,Z^*)$. 
Indeed, from (3.5) and the definition of
$\Lambda_r(y,K,t)$ it follows that $g_{t}{\bold A}_1(y)y\subset K$ for
all $y\in K\cup \{x\}$. Using  (3.6)  one can then inductively prove  that
$g_{jt}{\bold A}_j(y)y\subset K$ for all $y\in K\cup \{x\}$ and
$j\in\Bbb N$. This implies that 
$$
g_{jt}\ay(x)x\subset K\text{ for all }j \in\bn\,.\tag 3.7
$$
It remains to define  the set 
$
C = \bigcup_{s=-t}^{0} g_s K
$, which  is compact and disjoint from  $Z$ due to the
$F$-invariance of the latter. From
(3.7) it easily follows that  for any $h\in \ay(x)$, the orbit $Fhx$
is contained in $C$, and therefore is bounded and disjoint from $Z$.
\qed\enddemo

\heading{\S 4. \da\ and orbits of lattices}\endheading

\subhead{4.1}\endsubhead We return to the notation introduced in \S 1,
i.e.~put $G = SL_\mn(\br)\ltimes\br^\mn$, $\Gamma =
SL_\mn(\bz)\ltimes\bz^\mn$, $G_0 = SL_\mn(\br)$, $\Gamma_0 =
SL_\mn(\bz)$, $\Omega = \ggm$, the space of free lattices in
$\br^\mn$, $\la = \left(\matrix 
I_m & A  \\ 
0 & I_n
\endmatrix \right)$, $\lta = \langle \la,(\va,0)\t\rangle $, $
X = \text{diag}(\underbrace{\tfrac1m,\dots,\tfrac1m}_{\text{$m$ times}}, \underbrace{-\tfrac1n,\dots,-\tfrac1n}_{\text{$n$ times}})
$, 
$g_t = \exp(tX)$ and $F = \{g_t\mid t\ge 0\}$. 

As is mentioned in the introduction to \cite{KM1},
$\{\la\mid\text{\amr}\}$ is the $F$-\ehs\ of $G_0$. Similarly, one has 

\proclaim{Lemma}  $\{\lta\mid\text{\amtr}\}$ is the $F$-\ehs\ of $G$. \endproclaim

\demo{Proof} It is a straightforward computation that $\text{ad}\,X$ sends 
$\left\langle\left(\matrix
B & A  \\
C & D
\endmatrix \right),\left(\matrix
\va  \\
\vc
\endmatrix \right)\right\rangle\in\g$ (here $A\in M_{m,n}(\br)$, $B\in M_{m,m}(\br)$, $C\in M_{n,m}(\br)$, $D\in M_{n,n}(\br)$, $\va\in\br^m$, $\vc\in\br^n$) to the element 
$$
\left\langle\left(\matrix
0 & (\frac1m + \frac1n)A  \\
-(\frac1m + \frac1n)C & 0
\endmatrix \right),\left(\matrix
\frac1m \va  \\
-\frac1n\vc
\endmatrix \right)\right\rangle\,.\  \qed
$$
\enddemo

\subhead{4.2}\endsubhead We also need 
 
\proclaim{Lemma}  The action of $F$ on $\Omega$ is mixing. \endproclaim

\demo{Proof} Since $\br^\mn$ is the only nontrivial closed normal
subgroup of $G$, the \hs\ $\Omega$ has no nontrivial Euclidean
quotients  and its 
maximal semisimple quotient is equal to $G/\Delta$, where $\Delta =
\Gamma_0\ltimes\br^\mn$. Denote by  $p$ the quotient map $G\to
G_0\cong G/\br^\mn$. Then $G/\Delta$, as a $G$-space, is
$p$-equivariantly isomorphic to $G_0/\Gamma_0$,  and
clearly  $p(F)$  is not relatively compact in $G_0$. It
follows from  Moore's theorem that the  $F$-action  on $G/\Delta$ is
mixing, therefore,  by
Proposition 2.2, so is  the  $F$-action  on $\Omega$. \qed 
\enddemo

\subhead{4.3}\endsubhead We are now going to connect \di\ properties
of \amtr\ with orbit properties of $\lta\bz^\mn$. 
For comparison, let us first state Dani's
correspondence \cite{D1, Theorem 2.20} for homogeneous approximation. 

\proclaim{Theorem}  \amr\ is \ba\ iff there exists $\ve > 0$ such that $\|g_t\la\vv\| \ge \ve$ for all $t\ge 0$ and $\vv\in\bz^\mn\nz$.  \endproclaim

In view of Mahler's compactness criterion, the latter assertion is
equivalent to the orbit $F\la\bz^\mn$ being bounded (in other words,
to $\la\bz^\mn$ being an element of $E(F,\{\infty\})$). Therefore, as
is mentioned in \cite{K3},  one can use the result of \cite{KM1} to get
an alternative proof of Schmidt's theorem on thickness of the
set of  
\ba\ systems of linear forms. In order to move to affine forms, we
need an inhomogeneous analogue of the above criterion:

\proclaim{4.4. Theorem}  \amtr\ is  irrational  and  \ba\ iff 
$$
\exists\,\ve > 0\text{ such that }\|g_t\lta\vv\| \ge \ve\text{ for all
}t\ge 0\text{  and }\vv\in\bz^\mn\,.\tag 4.1 
$$
  \endproclaim

\demo{Proof} We essentially follow the argument of \cite{K1, Proof of
Proposition 5.2(a)}. Write $\vv = (\vp,\vq)\t$, where $\vp\in\bz^m$ and $\vq\in\bz^n$. Then 
$$
g_t\lta\vv = \big(e^{t/m}(A\vq + \va + \vp),e^{-t/n}\vq\big)\t\,.
$$
 This shows that (4.1) does not hold iff  there exist sequences $t_j
\ge 0$, $\vp_j\in\bz^m$ and $\vq_j\in\bz^n$ such that 
$$ 
e^{t_j/m}( A\vq_j+ \va + \vp_j)\to 0 
\tag 4.2a
$$
and
$$ 
e^{-t_j/n}\vq_j  \to 0\tag 4.2b
$$
 as  $j\to\infty$. On the other hand,  $\aa$  is \wa\ iff there exist sequences  $\vp_j\in\bz^m$ and $\vq_j\in\bz^n$, $j\in\bn$,  such that $\vq_j\to\infty$ and 
$$
\| A\vq_j + \va + \vp_j\|^m\|\vq_j\|^n\to 0\text{ as } j\to\infty\,.\tag 4.3
$$

We need to prove that (4.1) does not hold if and only if $\aa$  is either
rational or \wa. If  $\aa$  is rational, one can take $\vp_j = \vp_0$
and $\vq_j = 
\vq_0$, with $\vp_0$ and $\vq_0$ as in (1.2), and  arbitrary $t_j
\to \infty$; then the left hand side of (4.2a) is zero, and (4.2b) is
satisfied as well.  On the other hand, for   irrational and \wa\ $\aa$
one can define 
$e^{t_j}\df \sqrt{\|\vq_j\|^n / \| A\vq_j + \va + \vp_j\|^m}$ and
check that  
 (4.2ab) holds. 

Conversely, multiplying the norm of the left hand side of (4.2a) risen
to the $m$th power and  the norm of the left hand side of (4.2b)
risen to the $n$th power, one  immediately sees that (4.3) follows
from (4.2ab). It remains to observe that  either $\aa$  is 
rational or  $A\vq_j + \va + \vp_j$ is never zero, therefore (4.2a)
forces the sequence 
$\vq_j$ to tend to infinity.\qed
\enddemo

\subhead{4.5}\endsubhead Recall that we denoted by $\Omega_0$ the 
set of ``true'' (containing the zero vector) lattices in $\br^\mn$. 
It is now easy to complete the

\demo{Proof of Theorem 1.7} Take a \wa\ \amtr\  such that
$\lta\bz^\mn$ belongs to $E\big(F,(\Omega_0)^*\big)$. In view of the above criterion, 
$$
\exists \text{ a sequence }\Lambda_j\in F\lta\bz^\mn\text{ and vectors
}\vv_j\in\Lambda_j\text{ with }\vv_j\to 0\,.\tag 4.4
$$
 Since $F\lta\bz^\mn$ is relatively compact, one can without loss of generality assume that there exists $\Lambda\in\Omega$ with  $\Lambda_j\to\Lambda$ in the topology of $\Omega$. Clearly the presence of arbitrarily small vectors in the lattices $\Lambda_j$ forces $\Lambda$ to contain $0$, i.e.~belong to $\Omega_0$, which is a contradiction. 
\qed
\enddemo
 
\example{4.6. Remark} Note that the converse to Theorem 1.7 is not
true: by virtue of Theorem 4.4, any rational \amtr\ satisfies (4.4),
hence $\lta\bz^\mn$ is not in $E\big(F,(\Omega_0)^*\big)$. Restriction to the
irrational case gives a partial converse: indeed, the above proof
basically shows that the existence of a limit point $\Lambda\in \Omega_0$ of
the orbit $F\lta\bz^\mn$ violates (4.1); hence $\lta\bz^\mn$ belongs
to $E(F,\Omega_0)$ whenever $\aa$ is  irrational  and  \ba. But the orbit
$F\lta\bz^\mn$  
does not have to be bounded, as can be shown using the explicit
construction given by Kronecker's Theorem (see Example 1.4). Perhaps the
simplest possible example is the  irrational \ba\ form $\aa = \langle
0,1/2\rangle$ (here $m = 
n = 1$):  it is easy to see that the orbit 
$$
\text{diag}(e^t,e^{-t})\lta\bz^2
= \big\{\big(e^t(p + 1/2),e^{-t}q\big)\t\mid p,q\in\bz\big\}
$$
has no limit points in the space of free lattices in $\br^2$. 
\endexample


 \subhead{4.7}\endsubhead We conclude this section with the

\demo{Proof of Theorem 1.5} Observe that  $\Omega_0 = G_0\bz^\mn$ is
the orbit of a
proper subgroup  of $G$  containing $F$, which makes  it
null and $F$-invariant subset of $\Omega$. The fact that $\Omega_0$ is
closed is also straightforward. From Theorem
1.6 and Lemmas 4.1 and 4.2   it follows that the set
$\{\text{\amtr}\mid \lta\bz^\mn\in E\big(F,(\Omega_0)^*\big)\}$ is
thick in $\mtr$. In 
view of Theorem 1.7, systems of forms which belong to the latter set
are \ba, hence the thickness of the set $\batmn$. \qed 
\enddemo
 

\heading{\S 5. Concluding remarks and open questions}\endheading

It is worthwhile to look at the main result of this paper in the
context of other results in inhomogeneous \da. In what follows,
$\psi:\bn\mapsto(0,\infty)$ will be a
non-increasing function, and we will say, following \cite{KM3},
that a system of affine forms given by \amtr\ is {\it
$\psi$-approximable\footnote{We are grateful to M.$\,$M.$\,$Dodson for
 a permission to modify his terminology introduced in
\cite{Do1}. In our opinion, the use of (5.1) instead of traditional $
\|A\vq + \va + \vp\|  \le \psi(\|\vq\|)$ makes the structure,
including the connection with homogeneous flows, more transparent. See
\cite{KM2,KM3} for justification.}\/} if  
there exist infinitely many $\vq\in \bz^n$ such that  
$$
\|A\vq + \va + \vp\|^m  \le \psi(\|\vq\|^{n})\tag 5.1
$$
for some $\vp\in\bz^m$. Denote  by $\Cal W_{m,n}(\psi)$ the set of
$\psi$-approximable systems \amtr. 

The main result in the present paper is {\it doubly metric\/} in its
nature; that 
is, the object of study is the set of {\it all\/} pairs $\aa$. On
the contrary, in {\it singly metric\/} inhomogeneous problems one is
interested in the set of pairs $\aa$ where $\va$ is fixed. For
example, the doubly metric inhomogeneous Khintchine-Groshev Theorem
\cite{C, Chapter 
VII, Theorem II} says that  
$$
\text{the set }\Cal W_{m,n}(\psi)\text{ has}\quad\cases \text{full measure
\ if} \quad\sum_{k = 1}^\infty {\psi(k)} = \infty \\
\text{zero measure
if} \quad\sum_{k = 1}^\infty {\psi(k)} < \infty\,.
\endcases
$$
The
singly metric strengthening is due to Schmidt. For
$\va\in\br^m$, 
denote by 
$\Cal W_{m,n}(\psi, \va)$ the set of matrices \amr\ such that $\aa \in\Cal
W_{m,n}(\psi)$. It
follows from the main result of  \cite{S1} that given any
$\va\in\br^m$,  
$$
\text{the set }\Cal
W_{m,n}(\psi,\va)\text{ has}\quad\cases \text{full measure
\ if} \quad\sum_{k = 1}^\infty {\psi(k)} = \infty \\
\text{zero measure
if} \quad\sum_{k = 1}^\infty {\psi(k)} < \infty\,.
\endcases
$$
 Denote $\psi_0(x) = 1/x$. A quick comparison of (5.1) with (1.1)
shows that $\aa$
is \ba\ iff it is not $\ve\psi_0$-approximable for some $\ve >
0$; in other words, $\batmn = \mtr \smallsetminus \bigcup_{\ve >
0}\Cal W_{m,n}(\ve\psi_0)$. Therefore it follows from the above result
that for any  
$\va\in\br^m$,  the set
$$
\bamn(\va)\df \{\text{\amr}\mid \aa\in \batmn \}
$$
 has zero measure. 

Another class of related problems involves a connection between the
rate 
of decay of 
$\Cal W_{m,n}(\psi)$. The corresponding result in homogeneous
approximation is called the Jarnik-Besicovitch Theorem, see
\cite{Do1}. The doubly metric inhomogeneous version was done by
M.$\,$M.$\,$Dodson \cite{Do2} and H.$\,$Dickinson \cite{Di}, and recently
J.$\,$Levesley \cite{L} obtained  a singly metric strengthening. More
precisely he proved that for any
$\va\in\br^m$,  
$$
\text{dim}\big(\Cal
W_{m,n}(\psi,\va)\big) = \cases mn\big(1 - \frac{\lambda - 1}{m
+ n\lambda}\big)\text{ if } \lambda > 1 \\ mn
\qquad\qquad\quad\,\text{ if } \lambda \le 1\,,
\endcases
$$
where $\lambda = \liminf_{k\to\infty}
\left({\log\big(1/\psi(k)\big)}/{\log k}\right)$ is the {\it lower order\/} of the
function $1/\psi$. 

In view of the aforementioned results, one can ask whether it is
possible to prove that $\bamn(\va)$ is
thick in $\mr$ for every $\va\in\br^m$. (It is a consequence of Theorem
1.7 and slicing  
properties of the \hd\ 
that vectors $\va\in\br^m$ such that $\bamn(\va)$ is
thick form a thick subset
of $\br^m$.)

Note also that in the paper \cite{S3},  Schmidt proved
that  $\bamn$   is a {\it winning\/} (a property stronger than
thickness, cf.~\cite{S2, D2}) subset of $\mr$.  It is not clear to the
author whether Schmidt's methods can be modified to allow treatment of
inhomogeneous problems. It seems natural to conjecture that   $\batmn$
is a  winning subset of $\mtr$, and, moreover, that $\bamn(\va)$ is a
winning subset of $\mr$  for every $\va\in\br^m$.   This seems to be an
interesting and  challenging 
problem in 
metric number theory. 


\heading{  Acknowledgements}\endheading
The author is grateful to Professor G.$\,$A.$\,$Margulis for his interest
in 
this problem, and to the referee for useful comments and suggestions. 



\enddocument